\documentclass[11pt]{article}
\usepackage{amsmath}
\usepackage{amsfonts}
\usepackage{amsthm}
\usepackage{pictex}
\usepackage{graphicx}

\def\diag{\mathop{\mathrm{diag}}\nolimits}

\newtheorem{theorem}{Theorem}[section]
\newtheorem{lemma}[theorem]{Lemma}
\newtheorem{Remark}[theorem]{Remark}
\newenvironment{remark}{\begin{Remark}\rm}{\end{Remark}}
\numberwithin{equation}{section}

\begin{document}
\title{A Riemann-Hilbert problem for biorthogonal polynomials}
\author{A. B. J. Kuijlaars\thanks{Supported  by FWO research projects
G.0176.02 and G.0455.04.} \ \ and \
K. T-R McLaughlin\thanks{Supported by NSF grants \#DMS-9970328 and \#DMS-0200749.}}
\maketitle

\begin{abstract}
We characterize the biorthogonal polynomials that appear in
the theory of coupled random matrices via a Riemann-Hilbert
problem. Our Riemann-Hilbert problem is different from the
ones that were proposed recently by Ercolani and McLaughlin,
Kapaev, and Bertola et al. We believe that our formulation
may be tractable to asymptotic analysis.
\end{abstract}

\section{Introduction}
The biorthogonal polynomials that appear in the theory of coupled random
matrices \cite{BEH1,EM,EyMeh,Mehta2} are characterized by the property that
\begin{equation} \label{biorthogonality}
    \iint p_k(x) q_j(y) e^{-V(x)-W(y)+2\tau xy} dx dy = 0,
    \qquad \mbox{ if } j \neq k,
\end{equation}
where $p_k$ and $q_j$ are polynomials of exact degrees $k$ and
$j$, respectively. In (\ref{biorthogonality}) we have that $V, W :
\mathbb R \to \mathbb R$ are given functions with sufficient
increase at infinity so that the integrals converge, and $\tau \neq 0$
is a nonzero coupling constant. The integration is over $\mathbb R^2$.

Ercolani and McLaughlin \cite{EM} showed that the two sequences of
biorthogonal polynomials $(p_k)$ and $(q_j)$ exist, that they
are unique, and moreover, that $p_k$ has exactly $k$ simple
real zeros, see also \cite{Xu}.
They also gave a Riemann-Hilbert formulation for the biorthogonal
polynomials which is non-local in character. Recently, for the case that $V$
and $W$ are polynomials, Kapaev \cite{Kapaev}
and Bertola et al.\ \cite{BEH2}  gave local Riemann-Hilbert
problems. If $d = \deg W$, then the
Riemann-Hilbert problems for $p_k$  are formulated
for $d \times d$-matrix valued functions in \cite{BEH2, Kapaev}.

In this note we derive a different Riemann-Hilbert problem.
Our Riemann-Hilbert problem is based on the fact that the biorthogonal polynomials
can be characterized as multiple orthogonal polynomials (see below).
The formulation of a Riemann-Hilbert problem for multiple orthogonal
polynomials is due to Van Assche et al. \cite{VAGK}

An outstanding problem in random matrix theory is to provide a rigorous asymptotic
of eigenvalue statistics for coupled random matrices. The basic example is
the so-called $2$-matrix model in which we have
a probability measure on pairs $(M_1, M_2)$ of Hermitian $N \times N$ matrices
of the form
\[ \frac{1}{Z_N} \exp(- Tr(V(M_1) + W(M_2) - 2\tau M_1 M_2)) dM_1 dM_2. \]
Statistical quantities on eigenvalues of $M_1$ and $M_2$ can be expressed
in terms of the biorthogonal polynomials $p_k$ and $q_j$ given by (\ref{biorthogonality}),
see \cite{Mehta2, EyMeh}.
The connection to biorthogonal polynomials would be very useful, if one has,
in addition, a complete asymptotic description of the biorthogonal polynomials.
Then it would be possible to compute eigenvalue statistices in the large $N$ limit.
Indeed, although the calculations are somewhat involved, this has been carried out
in the Gaussian case $V(x) = x^2$, $W(y) = ay^2$, in \cite{EM}. Rigorous asymptotics
for biorthogonal polynomials with more general functions $V$ and $W$ are not
known.

In the $1$-matrix case, the statistical quantities on eigenvalues are given
in terms of orthogonal polynomials \cite{Mehta}, which have been characterized
by a Riemann-Hilbert problem \cite{FIK}. The steepest descent / stationary
phase method for Riemann-Hilbert problems was applied with great
success to orthogonal polynomials \cite{DKMVZ1,DKMVZ2,KMVV}.
As a result, the large $N$ asymptotics
of $1$-matrix models could be carried out in great detail, which in
particular provided a proof of the universality of eigenvalue spacings
for a large class of matrix models \cite{BI,DKMVZ1,KV}.
So there is hope that a similar asymptotic analysis of a Riemann-Hilbert
problem for biorthogonal polynomials will lead to large $N$ asymptotics
for $2$-matrix models. The formulation of a suitable Riemann-Hilbert problem is
only a first step in this direction.

For simplicity and clarity we formulate and prove our Riemann-Hilbert problem
for the biorthogonal polynomial $p_k$ for the first non-trivial case.
This is the case where $W$ is a
polynomial of degree $4$. Indeed, if $W$
is a polynomial of degree $2$, say $W(y) = y^2 + 2by + c$,
then $p_k$ is the orthogonal polynomial with respect to the weight
$e^{-V(x) + (\tau x-b)^2}$ on $\mathbb R$ and so there is
a Riemann-Hilbert problem for $p_k$ \cite{Deift,FIK}.

Thus we assume that $W$ is a polynomial of degree $4$ and we define
\begin{equation} \label{wj(x)}
    w_j(x) = \int y^j e^{-V(x)-W(y)+2\tau xy} dy, \qquad j =0,1,2.
\end{equation}
These functions appear in the formulation of our Riemann-Hilbert
problem.

\subsection*{Riemann-Hilbert problem for $Y$}
The problem is to find a $4\times 4$ matrix valued function
$Y : \mathbb C \setminus \mathbb R \to \mathbb C^{4 \times 4}$
having the following three properties.
\begin{enumerate}
\item[(a)] $Y$ is analytic on $\mathbb C \setminus \mathbb R$.
\item[(b)] $Y$ has boundary values on $\mathbb R$, denoted by $Y_+$
and $Y_-$, so that
\begin{equation} \label{jumpY}
    Y_+(x) = Y_-(x) \begin{pmatrix} 1 & w_0(x) & w_1(x) & w_2(x) \\
    0 & 1 & 0 & 0 \\
    0 & 0 & 1 & 0 \\
    0 & 0 & 0 & 1 \end{pmatrix}, \qquad x \in \mathbb R.
\end{equation}
\item[(c)] As $z \to \infty$, we have
\begin{equation} \label{asympY}
    Y(z) = \left(I + O\left(\frac{1}{z}\right)\right)
    \begin{pmatrix} z^k & 0 & 0 & 0 \\
    0 & z^{-n_0} & 0 & 0 \\
    0 & 0 & z^{-n_1} & 0 \\
    0 & 0 & 0 & z^{-n_2} \end{pmatrix},
\end{equation}
where $k \in \mathbb N_0$, $n_0 = \left[\frac{k+2}{3}\right]$,
$n_1 = \left[\frac{k+1}{3}\right]$, and $n_2 = \left[\frac{k}{3} \right]$.
Here $[\cdot]$ denotes the integer part.
(Note that $k = n_0 + n_1 + n_2$.)
\end{enumerate}

The main result of this paper is that the Riemann-Hilbert
problem for $Y$ has a unique solution and that its $(1,1)$ entry $Y_{11}$
is equal to the monic biorthogonal polynomial $p_k$.
In what follows we use $C(f)$ defined by
\[ C(f) (z) = \frac{1}{2\pi i} \int \frac{f(x)}{x-z} dx,
    \qquad z \in \mathbb C \setminus \mathbb R, \]
to denote the Cauchy transform of a function
$f: \mathbb R \to \mathbb R$.

\begin{theorem}
Let the functions $w_j$, $j=0,1,2$, be given by {\rm(\ref{wj(x)})} and
let $k \in \mathbb N_0$.
Then the above Riemann-Hilbert problem for $Y$ has a unique solution given
by
\begin{equation} \label{formulaY}
    Y = \begin{pmatrix}
    p_k & C(p_k w_0) & C(p_k w_1) & C(p_k w_2) \\[10pt]
    p_{k-1}^{(0)} & C(p_{k-1}^{(0)} w_0) &  C(p_{k-1}^{(0)} w_1) &
      C(p_{k-1}^{(0)} w_2)\\[10pt]
    p_{k-1}^{(1)} & C(p_{k-1}^{(1)} w_0) &  C(p_{k-1}^{(1)} w_1) &
      C(p_{k-1}^{(1)} w_2)\\[10pt]
    p_{k-1}^{(2)} & C(p_{k-1}^{(2)} w_0) &  C(p_{k-1}^{(2)} w_1) &
      C(p_{k-1}^{(2)} w_2)
      \end{pmatrix}
\end{equation}
where $p_k$ is the monic  polynomial of degree $k$ satisfying
{\rm(\ref{biorthogonality})}, and $p_{k-1}^{(j)}$, $j=0,1,2$, are three
polynomials of degrees $\leq k-1$.
\end{theorem}

\begin{remark}
There is an immediate extension to polynomials $W$ of arbitrary degree.
If $W$ is a polynomial of degree $d$, then the Riemann-Hilbert problem
is for $Y : \mathbb C \setminus \mathbb R \to \mathbb C^{d\times d}$
so that
\begin{equation} \label{RHgen1}
    Y \mbox{ is analytic on $\mathbb C \setminus \mathbb R$}.
\end{equation}
The jump condition uses the $d-1$ functions
$w_j(x) = \int y^j e^{-V(x)-W(y)+2\tau xy} dy$, $j=0,1,\ldots,d-2$,
and is given by
\begin{equation} \label{RHgen2}
    Y_+(x) = Y_-(x)
    \begin{pmatrix}
    1 & w_0(x) & w_1(x) & \cdots & \cdots & w_{d-3}(x) & w_{d-2}(x) \\
    0 & 1      & 0      & \cdots & \cdots & 0          & 0 \\
    0 & 0      & 1      & 0      &        &            & 0 \\
    \vdots&\vdots& 0    & \ddots & \ddots &            & \vdots \\
    \vdots&\vdots&      & \ddots & \ddots & \ddots     & \vdots \\
    0     & 0    &      &        & \ddots & \ddots     & 0 \\
    0     & 0    & 0    & \cdots & \cdots & 0          & 1
    \end{pmatrix},
\end{equation}
for $x \in \mathbb R$.
The asymptotic condition is
\begin{equation} \label{RHgen3}
    Y(z) = \left(I + O\left(\frac{1}{z}\right)\right)
    \diag\left(z^k, z^{-n_0}, z^{-n_1}, \ldots, z^{-n_{d-2}}\right)
\end{equation}
where $n_j = \left[\frac{k+d-2-j}{d-1}\right]$ for $j=0,1,\ldots,d-2$,
and $\diag\left( \cdot\right)$ denotes a diagonal matrix.

The Riemann-Hilbert problem (\ref{RHgen1})--(\ref{RHgen3})
has a unique solution and $Y_{11} = p_k$, where $p_k$
is the monic biorthogonal polynomial of degree $k$. The proof of the
general case follows along the same lines as the proof of the case
$\deg W = 4$ that we will present in Sections 2 and 3 below.
\end{remark}

\begin{remark}
By symmetry, there is a similar Riemann-Hilbert problem that characterizes
the other biorthogonal polynomial $q_j$, in the case that $V$ is a polynomial.
\end{remark}

\section{Multiple orthogonality}

We assume that $W$ is a polynomial of degree $4$. In this
section we will characterize the monic biorthogonal polynomial
$p_k$ of degree $k$ through a set of orthogonality relations
with respect to the three functions $w_j$, $j=0,1,2$.
As in Theorem 1.1, we use $n_0 = \left[\frac{k+2}{3}\right]$,
$n_1 = \left[\frac{k+1}{3}\right]$, and $n_2 = \left[\frac{k}{3} \right]$.
\begin{lemma}
We have
\begin{equation} \label{multipleorthogonality}
    \int p_k(x) x^i w_j(x) dx = 0,
    \qquad \mbox{for } i = 0, 1, \ldots, n_j-1, \
        j = 0,1,2,
\end{equation}
and these relations characterize the biorthogonal
polynomial $p_k$ among all monic polynomials of degree $k$.
\end{lemma}
\begin{proof}
Since $W$ is a polynomial of degree $4$, it is easy to see
that
\[ \frac{d^i}{dy^i} \left( y^j e^{-W(y)} \right) = \pi_{3i +j}(y) e^{-W(y)} \]
where $\pi_{3i+j}$ is a polynomial of exact degree $3i+j$.
For any function $f$, we then have if we integrate
by parts $i$ times
\begin{eqnarray} \nonumber
  \lefteqn{\iint f(x) \pi_{3i+j}(y) e^{-V(x)-W(y)+2\tau xy} dx dy \qquad \qquad} \\
  & & = \nonumber
    \int f(x) e^{-V(x)} \int \frac{d^i}{dy^i} \left(y^j e^{-W(y)} \right) e^{2\tau xy} dy dx \\
    & &  = \nonumber
    (-1)^i \int f(x) e^{-V(x)} \int y^j e^{-W(y)} \frac{d^i}{dy^i}\left(e^{2\tau xy}\right) dy dx \\
    & & = \nonumber
    (-2\tau)^i \int f(x) x^i e^{-V(x)} \int y^j e^{-W(y)+2\tau xy} dy dx \\
    & & = \label{integrationbyparts}
    (-2\tau)^i \int f(x) x^i w_j(x) dx.
\end{eqnarray}

If $f$ is the biorthogonal polynomial $p_k$ then the left-hand side of
(\ref{integrationbyparts}) is zero if $3i + j < k$. This corresponds exactly
with $i \leq n_j-1$ for $j=0,1,2$, so that by the right-hand side
we have the relations (\ref{multipleorthogonality}).

Conversely, if $f$ is a monic polynomial of degree $k$ that
satisfies the relations
$\int f(x) x^i w_j(x) dx = 0$ for $i=0,\ldots,n_j-1$, $j=0,1,2$,
then the left-hand side of (\ref{integrationbyparts})
is zero for $i \leq n_j-1$ and $j=0,1,2$.
The polynomials $\pi_{3i+j}$ with $i \leq n_j-1$ and $j=0,1,2$
are a basis for the polynomials of degree $\leq k-1$. Hence
$f$ is the biorthogonal polynomial $p_k$.
\end{proof}

\begin{remark}
The relations (\ref{multipleorthogonality}) are called multiple orthogonality
relations of type II, see \cite{Aptekarev,ABV,NS,VAC,VAGK} for more on this subject.
\end{remark}
\section{Proof of Theorem 1.1}

\begin{proof}
We first establish uniqueness  in the standard way.
The only thing we have to observe is that both the jump matrix in (\ref{jumpY})
and the diagonal matrix in the right-hand side of (\ref{asympY})
have determinant one. Then the proof of uniqueness follows
as in \cite[Section 3.2]{Deift}.

We now prove that $Y$ given by (\ref{formulaY}) satisfies
the Riemann-Hilbert problem.
First we consider the first row of $Y$.
The conditions (\ref{jumpY}) and (\ref{asympY}) give for the
$(1,1)$ entry
\[ Y_{11,+} = Y_{11,-}, \quad \mbox{ and } \quad Y_{11}(z) = z^k + O(z^{k-1}), \]
These conditions are clearly satisfied if $Y_{11} = p_k$, since $p_k$
is a monic polynomial of degree $k$.

For the other entries in the first row, the jump condition (\ref{jumpY})
then is
\[ Y_{1j,+} = Y_{1j,-} + Y_{11,-}w_{j-2} = Y_{1j,-} + p_k w_{j-2} \qquad j=2,3,4.  \]
By the Sokhotskii-Plemelj formula, this is satisfied
by $Y_{1j} = C(p_k w_{j-2})$.
The asymptotic condition (\ref{asympY}) is
\begin{equation} \label{asympY1j}
    Y_{1j}(z) = O(z^{-n_{j-2}-1}) \qquad \mbox{ as } z \to \infty
    \end{equation}
and we have to check that this is satisfied for $Y_{1j} = C(p_kw_{j-2})$.
If we use
\[ \frac{1}{x-z} = - \sum_{i=0}^{n_{j-2}-1} \frac{x^i}{z^{i+1}}
    + \frac{x^{n_{j-2}}}{z^{n_{j-2}}} \frac{1}{x-z} \]
then we see that for $j=2,3,4$,
\begin{eqnarray} \nonumber
    C(p_k w_{j-2})(z) & = & \frac{1}{2\pi i} \int \frac{p_k(x) w_{j-2}(x)}{x-z} \\
    & = & \nonumber - \sum_{i=0}^{n_{j-2}-1} \left(\frac{1}{2\pi i} \int p_k(x) x^i w_{j-2}(x) dx\right) z^{-i-1} \\
    & & \
    + \left( \frac{1}{2\pi i} \int \frac{p_k(x) w_{j-2}(x) x^{n_{j-2}}}{x-z} dx \right)
    z^{-n_{j-2}}. \label{expansionCpkwj}
\end{eqnarray}
Because of the multiple orthogonal relations (\ref{multipleorthogonality}) satisfied
by $p_k$, (\ref{expansionCpkwj}) reduces to
\[ C(p_k w_{j-2})(z) = \left( \frac{1}{2\pi i} \int \frac{p_k(x) w_{j-2}(x) x^{n_{j-2}}}{x-z} dx \right)
    z^{-n_{j-2}} \]
which shows that (\ref{asympY1j}) is indeed satisfied if $Y_{1j} = C(p_kw_{j-2})$

Next we consider the second row of $Y$.
The conditions $Y_{21,+} = Y_{21,-}$ and $Y_{21}(z) = O(z^{k-1})$ as $z \to \infty$
are clearly satisfied if $Y_{21}$ is a polynomial $p_{k-1}^{(0)}$ of degree $\leq k-1$.
The jump conditions for the other entries in the second row
\[ Y_{2j,+} = Y_{2j,-} + Y_{1j,-}w_{j-2} = Y_{2j,-} + p_{k-1}^{(0)} w_{j-2} \]
are then also satisfied if $Y_{2j} = C(p_{k-1}^{(0)}w_{j-2})$ for $j=2,3,4$.
We need to be able to choose $p_{k-1}^{(0)}$ so that the asymptotic condition
(\ref{jumpY}) is also satisfied, which means that
\begin{equation} \label{asympCpk-1}
    \left\{ \begin{array}{ll}
    C(p_{k-1}^{(0)}w_0) = z^{-n_0} + O(z^{-n_0-1}), & \\[10pt]
    C(p_{k-1}^{(0)}w_j) = O(z^{-n_j-1}) & \mbox{ for } j=1,2.
    \end{array} \right.
\end{equation}
Expanding $C(p_{k-1}^{(0)}w_j)$ as in (\ref{expansionCpkwj}),
we see that (\ref{asympCpk-1}) is satisfied if $p_{k-1}^{(0)}$ is such that
\begin{eqnarray}
    \label{orthog1}
    \int p_{k-1}^{(0)}(x) x^i w_0(x) dx = &0, & \quad \mbox{for } i=0, \ldots, n_0-2, \\
    \label{orthog2}
    \int p_{k-1}^{(0)}(x) x^i w_0(x) dx = &- 2\pi i, & \quad \mbox{for } i = n_0-1, \\
    \label{orthog3}
    \int p_{k-1}^{(0)}(x) x^i w_j(x) dx = &0, & \quad \mbox{for } i=0, \ldots, n_j-1, j=1,2.
\end{eqnarray}
The conditions (\ref{orthog1}) and (\ref{orthog3}) give $k-1$ homogeneous
conditions on the $k$ free coefficients of $p_{k-1}^{(0)}$,
and so there exists a non-zero polynomial $p_{k-1}^{(0)}$
satisfying these conditions.
To be able to have (\ref{orthog2}) as well we must exclude the
possibility that
\begin{equation} \label{normalizepk-1}
    \int p_{k-1}^{(0)} (x) x^{n_0-1} w_0(x) dx = 0.
\end{equation}
However, if (\ref{normalizepk-1}) would hold, then
$p_k + p_{k-1}^{(0)}$ would be a monic polynomial of degree $k$
that satisfies the  multiple orthogonality relations (\ref{multipleorthogonality}),
which is impossible, since these relations characterize the biorthogonal
polynomial $p_k$ by Lemma 2.1. Thus (\ref{normalizepk-1}) cannot hold.
Then we can normalize $p_{k-1}^{(0)}$ by multiplying it with a suitable constant, so that
(\ref{orthog2}) is satisfied. This proves that we can indeed choose $p_{k-1}^{(0)}$
so that the second row of $Y$ satisfies all the conditions imposed by the Riemann-Hilbert
problem.

In exactly the same way, we handle the third and fourth rows of $Y$.

This completes the proof of Theorem 1.1.
\end{proof}

\section{Conclusion}
We have characterized the biorthogonal polynomials that appear in
the theory of coupled random matrices via a Riemann-Hilbert problem
which is differnt from the Riemann-Hilbert problems derived in
\cite{BEH2, Kapaev}.
Recent experience \cite{BK,KVW} with similar higher order Riemann-Hilbert
problems leads us to believe that our Riemann-Hilbert problem may be
tractable to asymptotic analysis. However, up to now we have not been
able to apply the steepest descent method succesfully to this problem,
and the actual asymptotic analysis of biorthogonal polynomials remains
a major open problem.

\obeylines
\texttt{
A. B. J. Kuijlaars
Department of Mathematics
Katholieke Universiteit Leuven
Celestijnenlaan 200B
B-3001 Leuven, BELGIUM,
E-mail: arno@wis.kuleuven.ac.be
\bigskip

K. T-R McLaughlin 
Department of Mathematics
University of North Carolina at Chapel Hill
Chapel Hill NC 27599, U.S.A.
E-mail: mcl@amath.unc.edu
}

\end{document}